\input amstex
\documentstyle{amsppt}
\document

\topmatter
\title Perverse Sheaves on Loop Grassmannians \\
 and Langlands Duality
\endtitle
\author Ivan Mirkovi\'c and Kari Vilonen
\endauthor
\address {Department of mathematics, University of Massachusetts, Amherst, MA
01002, USA}
\endaddress
\email mirkovic\@math.umass.edu
\endemail
\address{Department of mathematics, Brandeis University, Waltham, MA 02254,
USA}
\endaddress
\email vilonen\@math.brandeis.edu
\endemail
\thanks I. Mirkovi\'c was partially supported by NSF
\endgraf K.Vilonen was partially supported by NSA and NSF
\endthanks

\endtopmatter

\define\V{{\operatorname {Vec}_\Bbbk}}
\define\oh{{\operatorname H}}

\define\ps{\operatorname {P}_{\cs}}
\define\p {\operatorname {P}}
\define\pg#1{\operatorname P_{\GO}(\G,#1)}
\define\GO{{G(\Cal O)}}
\define\GK{{G(\Cal K)}}
\define\G{{\Cal G}}


\define\bc{{\Bbb C}}

\define\cs{{\Cal S}}
\define\ct{{\Cal T}}
\define\cf{{\Cal F}}
\define\ca{{\Cal A}}

\define\ck{{\Cal K}}
\define\co{{\Cal O}}

%


\define\ra{\rightarrow}


\define\bb{\underset}           

\define\barr{\overline}         
\predefine\ss{\S}               

\define\sub{\subseteq}          

\define\tim{\times}             






\define\12{ \frac{1}{2} }       
\define\cd{{\cdot}}             



\define\inv{{}^{-1}}

\define\Hom{\operatorname{Hom}}

\define\a1{{\A^1}}              


\redefine\AA{\Cal A}
\define\BB{\Cal B}

\define\GG{\Cal G}

\define\KK{\Cal K}

\define\OO{\Cal O}
\define\PP{\Cal P}

\define\k{{\Bbbk }}

\redefine\A{\Bbb A }

\redefine\C{\Bbb C}


\redefine\S{\Bbb S}



\define\fg{\frak g}




\define\la{\lambda }




















\define\PGG#1{\operatorname P_{\GO}(\GG,#1)}








\define\htt{\operatorname{ht}}


\subheading{\bf 1. Introduction}
\vskip .5cm

In this paper we outline a proof of a  geometric version of  the Satake
isomorphism. Namely, given a connected, complex algebraic reductive  group $G$
we show that the tensor category of representations of the dual group $\check
G$ is naturally equivalent to a certain category of perverse sheaves on the
loop Grassmannian of $G$. The tensor category structure on this category of
perverse sheaves is given  via a convolution product.

The above result is not new. It has been announced by Ginsburg in
\cite{G} and  some of the arguments in section 5 of  this paper are borrowed
from \cite{G}. However, at  crucial points our proof differs from Ginsburg's.
First, we use a more \lq\lq natural" commutativity constraint for the
convolution product. This commutativity constraint, explained in section 3,
is due to Drinfeld and was explained to us by Beilinson. Secondly, in section
4,  we give a direct geometric proof that the global cohomology functor is
exact and decompose this cohomology functor into a direct sum of weights
(Theorem 4.3). We completely avoid the use of the decomposition theorem of
\cite{BBD} which makes our techniques applicable to perverse sheaves with
coefficients over arbitrary commutative rings.

This note includes sketches of (some of the) proofs. The details, as well as
the generalization of the results from
$\Bbb C$-representations to representations over arbitrary fields and
commutative rings will appear elsewhere.

\vskip .5cm

\subheading{\bf 2. The Convolution Product}
\vskip .5cm

Let $G$ be a connected, complex algebraic reductive group. Denote by
$\Cal O = \Bbb C[[t]]$  the  ring of formal  power series in one variable and
by
$\Cal K = \Bbb C((t))$ its fraction field, the field of formal Laurent series.
The loop Grassmannian, as a set, is defined as $\G =
\GK/\GO$, where, as usual, $\GK$ and $\GO$ denote the sets of the
$\ck$-valued and the $\co$-valued points of $G$ respectively. The  sets
$\GK$,
$\GO$, and $\G$ have an algebraic structure as $\Bbb C$-spaces. The space
$\GO$ is a group scheme over $\Bbb C$ but  the spaces $\GK$ and $\G$ are only
ind-schemes\footnote{By an ind-scheme we mean an ind-scheme in a strict
sense, i.e., an inductive system of schemes where all maps are
closed embeddings.}. To see that
$\GK$ is an ind-scheme, one embeds $G$ in $SL_N(\bc)$. The filtration by order
of pole in $SL_N(\ck)$ induces a filtration of $\GK$ which exhibits
$\GK$ as an inductive limit of schemes. The filtration above is invariant under
the (right) action of $\GO$ on $\GK$ and thus, after taking the quotient of
$\GK$  by $\GO$ one gets a filtration of $\G$ which exhibits it as a union of
finite dimensional projective schemes. Furthermore, the morphism $\pi:\GK
\to
\G$ is locally trivial in the Zariski topology, i.e., there exists a Zariski
open subset $U\subset\G$ such that $\pi^{-1}(U)\cong U\times\GO$ and $\pi$
restricted to $ U\times\GO$ is simply projection to the first factor. For
details see for example \cite{BL1,LS}.

The group scheme $\GO$ acts on $\G$ with finite dimensional orbits. In order
to describe the orbit structure, let us fix a maximal torus $T\subset G$. We
write $W$ for the Weyl group and $X_*(T)$  for the coweights
$\operatorname{Hom}(\bc^*,T)$. Then the $\GO$-orbits on $\G$ are parametrized
by the $W$-orbits in $X_*(T)$, and given $\lambda\in X_*(T)$ the $\GO$-orbit
associated to it is $\G_\lambda = \GO
\cdot\lambda\subset\G$, where we have identified $X_*(T)$ as a subset of
$\GK$.

Let $\Bbbk$ be a field of characteristic zero, which we fix for the rest of the
paper. All sheaves that we encounter in this paper will be sheaves in the
classical topology.  We denote by $\pg \Bbbk$ the category of
$\GO$-equivariant perverse
$\Bbbk$-sheaves  on $\G$ with finite dimensional support and by
$\ps(\G,\Bbbk)$ the category of perverse $\Bbbk$-sheaves on $\G$ which are
constructible with respect to the orbit stratification $\cs$ of $\G$ and which
have finite dimensional support. We use the notational conventions of
\cite{BBD} for perverse sheaves, in particular, in order for the constant sheaf
on a
$\GO$-orbit $\G_\lambda$ to be perverse it has to be placed in degree $-\dim
\G_\lambda$.

\proclaim{Proposition 2.1} The forgetful functor $\pg \Bbbk \to
\ps(\G,\Bbbk)$ is an equivalence of categories.
\endproclaim

We will now put a tensor category structure on $\pg \Bbbk$ via the convolution
product. Consider the following diagram of maps (of sets)
$$
\G\times\G @<p<< \GK \times \G @>q>> \GK \times_\GO
\G @>m>> \G\,.
\tag2.2
$$  Here $\GK \times_\GO\G$ denotes the quotient of $\GK \times \G$ by
$\GO$ where the action is given on the $\GK$-factor via right multiplication
by an inverse and on the $\G$-factor by left multiplication. The  $p$ and $q$
are projection maps and $m$ is the multiplication map. All other terms in (2.2)
have been given a structure of an ind-scheme except $\GK \times_\GO\G$. The
description of this structure is easier in the global context of section 3
where it is a special case of a more general construction and thus we omit the
details here. We define the convolution product
$A_1 *A_2$ of
$A_1,A_2\in\pg \Bbbk$ by the formula
$$ A_1 *A_2 \ = \ Rm_*\tilde A \qquad\text{where \ $q^*\tilde A  =
p^*(A_1\boxtimes A_2)$}\,.
\tag2.3
$$ To make sense of this definition we first use the fact  that $p$   and
$q$ are locally trivial in the Zariski topology. This guarantees the existence
of
$\tilde A\in\p_\GO(\GK\times_\GO\G,\Bbbk)$. To see the local triviality of
$q$ one can use the same arguments as for example in
\cite{BL1,LS}, and as was pointed out above, the local triviality of $p$ is
proved in those references. It remains to show that
$Rm_*\tilde A \in\pg\Bbbk$. To that end we introduce the notion of a
stratified semi-small map.

Let us consider two complex stratified spaces $(Y,\ct)$ and $(X,\cs)$ and a
map $f:Y\to X$. We assume that the two stratifications are locally trivial with
connected strata and that $f$ is a stratified with respect to the
stratifications $\ct$ and $\cs$, i.e., that for any $T\in\ct$ the image
$f(T)$ is a union of strata in $\cs$ and for any $S\in\cs$ the map
$f|f^{-1}(S):  f^{-1}(S)
\to S$ is locally trivial in the stratified sense. We say that
$f$ is a stratified semi-small map if

$$
\aligned a) \ \ &\text{for any $T\in\ct$ the map $f|\bar T$ is proper}
\\ b) \ \ &\text{for any $T\in\ct$ and  any $S\in\cs$ such that $S\subset
f(\bar T)$ we have}
\\ &\dim(f^{-1}(x)\cap \bar T) \leq \frac 1 2 (\dim f(\bar T) - \dim S)
\\ &\text{for any (and thus all) $x\in S$\,. }
\endaligned
\tag2.4
$$ Next the notion of a small stratified map. We say that $f$ is a small
stratified map if there exists a (non-trivial) open stratified subset $W$ of
$Y$ such that
$$
\aligned a) \ \ &\text{for any $T\in\ct$ the map $f|\bar T$ is proper}
\\ b) \ \ &\text{the map $f|W:W\to f(W)$ is proper and has finite fibers}
\\ c) \ \ &\text{for any $T\in\ct$, $T\subset W$, and  any $S\in\cs$ such that
$S\subset f(\bar T)-f(T)$}
\\ &\text{we have}\ \ \dim(f^{-1}(x)\cap \bar T) \leq \frac 1 2 (\dim f(\bar T)
-
\dim S)
\\ &\text{for any (and thus all) $x\in S$\,. }
\endaligned
\tag2.5
$$

The result below follows directly from dimension counting:
\proclaim{Lemma 2.6} If $f$ is a semismall stratified map then
$Rf_*A\in\p_\cs(X,\Bbbk)$ for all $A\in\p_\ct(Y,\Bbbk)$\,. If $f$ is a small
stratified map then, with any $W$ as above, and any $A\in\p_\ct(W,\Bbbk)$, we
have $Rf_* j_{!*} A = \tilde j_{!*} f_*A$, where $j: W \hookrightarrow Y$ and
$\tilde j :f(W) \hookrightarrow X$ denote the two inclusions.

\endproclaim

We apply the above considerations, in the semismall case,  to our situation.
We take $Y=\GK\times_\GO\G$ and  choose $\ct$ to be the stratification whose
strata are $p^{-1}(\G_\lambda)\times_\GO\G_\mu$, for
$\lambda,\mu\in X_*(T)$\,. We also let $X=\G$, $\cs$ the stratification by
$\GO$-orbits, and choose
$f=m$. To conclude the constructions of the convolution product on $\pg
\Bbbk$ it suffices to  note that the sheaf $\tilde A$ is constructible with
respect to the stratification $\ct$ and appeal to the following

\proclaim{Theorem 2.7} The multiplication map
$\GK \times_\GO \G @>m>> \G$ is a stratified semi-small map with respect to
the stratifications above.
\endproclaim

For an outline of proof, see the appendix.

One can define the convolution product of three sheaves completely analogously
to (2.3). This gives an associativity constraint for the convolution product
thus giving $\pg \Bbbk$ the structure of an associative tensor category. In
the next section we construct a commutativity constraint for the convolution
product.

\vskip 1cm

\subheading{\bf 3. The Commutativity Constraint}
\vskip .5cm

In order to construct the commutativity constraint we will need to consider
the convolution product in the global situation. Let $X$ be a smooth curve
over the complex numbers. Let $x\in X$ be a closed point and denote by
$\Cal O_x$ the completion of the local ring at $x$ and by $\Cal K_x$ its
fraction field. Then the Grassmannian $\G_x = G(\Cal K_x)/G(\Cal O_x)$
represents the following functor from $\Bbb C$-algebras to sets :
$$ R \mapsto \{
\cf
\text{  a  $G$-torsor on $X_R$,
$\nu: G\times X^*_R \ra \cf|X^*_R$ a trivialization on $X^*_R$  }
\}\,.
\tag3.1
$$ Here the pairs $(\cf,\nu)$ are to be taken up to isomorphism, $X_R =
X\times\text{Spec}(R)$, and $X^*_R= (X-\{x\})\times \text{Spec}(R)$\,. For
details see for example \cite{BL1,BL2,LS}. We now globalize this construction
and at the same time form the Grassmannian at several points on the curve.
Denote the
$n$ fold product by
$X^n = X\times
\dots\times X$  and consider the functor
$$ R \mapsto \left\{\aligned &(x_1,\dots,x_n) \in X^n(R), \ \ \cf\text{ a
$G$-torsor on
$X_R$\,, }\\ &\text{$\nu_{(x_1,\dots,x_n)}$ a trivialization of
$\cf$ on $X_R - \cup {x_i}$}\endaligned\right\}\,.
\tag3.2
$$ Here we think of the points $x_i: \text{Spec}(R) \to X$ as subschemes of
$X_R$ by taking their graphs. One sees that the functor in (3.2) is
represented by an ind-scheme $\G_X^{(n)}$. Of course  $\G_X^{(n)}$ is an
ind-scheme over
$X^n$ and its fiber over the point
$(x_1,\dots,x_n)$ is simply $\prod_{i=1}^k \G_{y_i}$\,, where
$\{y_1,\dots,y_k\}=\{x_1,\dots,x_n\}$, with all the $y_i$ distinct. We write
$\G_X^{(1)} = \G_X$.

We will now extend the diagram of maps (2.2), which was used to define the
convolution product, to the global situation, i.e., to a diagram of ind-schemes
over $X^2$:
$$
\G_X\times\G_X @<p<< \widetilde{\G_X\times \G_X} @>q>> \G_X  \tilde
\times
\G_X  @>m>>\G_X^{(2)}\,.
\tag3.3
$$ Roughly, the diagram starts with a pair of torsors, each trivialized off
one  point. One chooses a trivialization of the first torsor near the second
point, and uses it to glue the torsors.

More precisely,
$\widetilde{\G_X\times \G_X}$ denotes the ind-scheme representing the functor
$$ R \mapsto \left\{\aligned &(x_1,x_2)\in X^2(R); \ \cf_1,\cf_2\text{
$G$-torsors on
$X_R$;\
$\nu_i$ a trivialization of }
\\ &\text{$\cf_i$ on $X_R -x_i$, for $i=1,2$; \ $\mu_1$ a trivialization of
$\cf_1$ on
$\widehat{(X_R)}_{x_2}$}
\endaligned
\right\},
\tag3.4
$$ where $\widehat{(X_R)}_{x_2}$ denotes the formal neighborhood of
$x_2$ in
$X_R$. The \lq\lq twisted product" $\G_X  \tilde \times \G_X $ is the
ind-scheme representing the functor
$$ R \mapsto\left\{\aligned &(x_1,x_2)\in X^2(R); \ \cf_1,\cf \text{
$G$-torsors on $X_R$; $\nu_1$ a trivialization }
\\ &\text{of
$\cf_1$ on $X_R-x_1$;} \ \eta : \cf_1 | (X_R -x_2) @>{\ \ \simeq \ \ }>>
\cf| (X_R -x_2)
\endaligned
\right\}\,.
\tag3.5
$$ It remains to describe the morphisms $p$, $q$, and $m$ in (3.3). Because all
the spaces in (3.3) are ind-schemes over $X^2$, and all the functors  involve
the choice of the same $(x_1,x_2)\in X^2(R)$ we omit it in the formulas below.
The morphism $p$ simply forgets the choice of $\mu_1$, the morphism
$q$ is given by the natural transformation
$$
 (\cf_1,\nu_1,\mu_1;\cf_2,\nu_2) \mapsto (\cf_1,\nu_1,\cf,\eta),
\tag3.6
$$ where $\cf$ is the $G$-torsor gotten by gluing $\cf_1$ on $X_R - x_2$ and
$\cf_2$ on $\widehat{(X_R)}_{x_2}$ using the isomorphism induced by
$\nu_2\circ\mu_1^{-1}$ between $\cf_1$ and $\cf_2$ on
$(X_R-x_2)\cap \widehat{(X_R)}_{x_2}$. The morphism
$m$ is given by the natural transformation
$$ (\cf_1,\nu_1,\cf,\eta) \mapsto (\cf,\nu)\,,
\tag3.7
$$ where $\nu = (\eta \circ \nu_1)|(X_R - x_1-x_2)$.

Next, the global analog of $\GO$ is the group-scheme $G_X^{(n)}(\Cal O)$ which
represents the functor
$$ R \mapsto \left\{\aligned &(x_1,\dots,x_n) \in X^n(R), \ \ \cf\text{ the
trivial
$G$-torsor on
$X_R$\,, }\\ &\text{$\mu_{(x_1,\dots,x_n)}$ a trivialization of
$\cf$ on $\widehat{(X_R)}_{(x_1\cup\dots\cup x_n)}$}\endaligned\right\}\,.
\tag3.8
$$

Just  as in section 2 we define the convolution product of $\BB_1,\BB_2\in
\PP_{G_X(\Cal O)}(\G_X,\Bbbk)$ by the formula
$$
\BB_1 \bb{X}\to* \BB_2 \ = \ Rm_*\tilde \BB \qquad
\text{where \ $q^*\tilde \BB = p^*(\BB_1\boxtimes \BB_2)$}\,.
\tag3.9
$$ Precisely as in section 2, the sheaf $\tilde \BB$ exists because $q$ is
locally, even in the Zariski topology, a product. Furthermore, the map $m$ is a
stratified small map -- regardless of the stratification on $X$. To see this,
let us denote by $\Delta \subset X^2$ the diagonal and set $U = X^2-\Delta$.
Then we can take, in  definition (2.5), as $W$ the locus of points lying over
$U$. That $m$ is small now follows as
$m$ is an isomorphism over $U$ and over points of $\Delta$ the map $m$
coincides with its analogue in section 2 which is semi-small by theorem 2.7.

Let us now, for simplicity, choose $X=\Bbb A^1$. Then the choice of a global
coordinate on $\Bbb A^1$, trivializes $\G_X$ over $X$; let us write $\rho :\G_X
\to \G$ for the projection. Let us denote $\rho^0 =
\rho^*[1] : \PP_\GO(\GG,\Bbbk) \to \PP_{G_X(\Cal O)}(\G_X,\Bbbk)$\,.  By
restricting
$\G_X^{(2)}$ to the diagonal $\Delta \cong X$ and to $U$, and observing that
these restrictions are isomorphic to $\G_X$ and to
$(\G_X\times\G_X)|U$ respectively, we get the following diagram
$$
\CD
\G_X @>i>> \G_X^{(2)} @<j<< (\G_X\times\G_X)|U
\\ @VVV @VVV @VVV
\\ X @>>> X^2 @<<<\ \ \ U\ \ \ .
\endCD
\tag3.10
$$

\proclaim{Lemma 3.11} For $\AA_1,\AA_2\in\PP_\GO(\GG,\Bbbk)$ we have
$$
\aligned &\text{a)}\qquad
\rho^0 \AA_1   \bb{X}\to * \rho^0 \AA_2 \ \cong \ j_{!*}\left( (
\rho^0\AA_1\boxtimes\rho^0\AA_2)|U \right)
\\ &\text{b)}\qquad\rho^0(\AA_1*\AA_2) \ \cong \  i^0 ( \rho^0 \AA_1
\bb{X}\to * \rho^0\AA_2 ) \,.
\endaligned
$$
\endproclaim Part a) of the lemma follows from smallness of $m$ and lemma 2.6.

Lemma 3.11 gives us the following sequence of isomorphisms:
$$
\gathered
\rho^0(\AA_1*\AA_2) \cong i^0j_{!*} \left(
(\rho^0\AA_1\boxtimes\rho^0\AA_2)|U \right)
\\
\cong  i^*j_{!*}((\rho^0\AA_2\boxtimes\rho^0\AA_1)|U)
\cong\rho^0(\AA_2*\AA_1)\,.
\endgathered
\tag3.12
$$ Specializing this isomorphism to (any) point on the diagonal yields a
functorial isomorphism between $\AA_1*\AA_2$ and $\AA_2*\AA_1$. We take this
isomorphism as our commutativity constraint.

\remark{Remark 3.13} The construction of the commutativity constraint can be
carried out in a more elegant way as follows. We first observe that the image
of the embedding $\rho^0 =\rho^*[1] : \PP_\GO(\GG,\Bbbk) \to
\PP_{G_X(\Cal O)}(\G_X,\Bbbk)$ consists precisely of objects in
$\PP_{G_X(\Cal O)}(\G_X,\Bbbk)$ which are \lq\lq constant" along $X$. This
subcategory of $\PP_{G_X(\Cal O)}(\G_X,\Bbbk)$ coincides with
$\PP_{\tilde G_X(\Cal O)}(\G_X,\Bbbk)$, where $\tilde G_X(\Cal O)$ denotes the
semi direct product of $G_X(\Cal O)$ and the groupoid which consists of pairs
of points $(x,y)\in X\times X$ together with an isomorphism between the formal
neighborhood of $x$ and the formal neighborhood of $y$. Now
$\rho^0 =\rho^*[1] : \PP_\GO(\GG,\Bbbk) \to\PP_{\tilde G_X(\Cal
O)}(\G_X,\Bbbk)$ is an equivalence whose inverse is $i^0 = i^*[-1]$, where
$i: \G_x \hookrightarrow \G_X$ is the inclusion. If $X$ is an arbitrary smooth
curve then the functor $i^0: \PP_{\tilde G_X(\Cal O)}(\G_X,\Bbbk) \to
\PP_\GO(\GG,\Bbbk)$ still has meaning and is an equivalence of categories. It
is clear that the convolution product (3.9) gives us a convolution product on
the category $\PP_{\tilde G_X(\Cal O)}(\G_X,\Bbbk)$. Thus, we can give the
construction of the commutativity constraint in terms of $\PP_{\tilde G_X(\Cal
O)}(\G_X,\Bbbk)$ and $i^0$ without specializing to $X= \Bbb A^1$ and choosing
a global coordinate.
\endremark
\vskip 1cm

\subheading{\bf 4. The Fiber Functor}
\vskip .5cm

Let $\V$ denote the category of finite dimensional vector spaces over
$\Bbbk$. Let us consider the global cohomology functor $\Bbb H^*:
\pg\Bbbk \to \V$, where we ignore the grading on global cohomology.

\proclaim{Proposition 4.1} The functor  $\Bbb H^*: \pg\Bbbk \to \V$ is a
tensor  functor.
\endproclaim

Let $r$ denote the map $r:\G_X^{(2)} \to X^2$. That $\Bbb H^*$ is tensor
functor follows immediately from
$$
\aligned &\text{a)}\ \ Rr_*(\rho^0(A_1)*_X\rho^0(A_2))|U  \ \ \text{is the
constant sheaf}\ \
\Bbb H^*(A_1)\otimes \Bbb H^*(A_2)\,.
\\ &\text{b)}\ \ Rr_*(\rho^0(A_1)*_X\rho^0(A_2))|\Delta = \rho^0(\Bbb
H^*(A_1*A_2))
\\ &\text{c)}\ \ Rr_*(\rho^0(A_1)*_X\rho^0(A_2)) \ \ \text{is a constant sheaf}
\endaligned
\tag4.2
$$ The claims a) and b) follow from lemma 3.11. It remains to note that, in the
notation of formula (3.9), the sheaf $R(r\circ m)_*\tilde B$ is constant; this
implies c).

We now come to the main technical result of this paper. In order to state it we
will fix some further notation. We choose a Borel subgroup $B\subset G$ which
contains the maximal torus $T$. This, of course, determines a choice of
positive roots. Let $N$ denote the unipotent radical of $B$. As usual, we
denote by
$\rho$ half the sum of positive roots of $G$. For any $\nu\in X_*(T)$ we write
$\htt(\nu)$ for the height of $\nu$ with respect to $\rho$. The
$N(\ck)$-orbits on
$\G$ are parametrized by $X_*(T)$; to each $\nu\in X_*(T) = \Hom(\Bbb C^*,T)$
we associate the $N(\ck)$-orbit $S_\nu =_{\text{def}}N(\ck)\cdot\nu$. Note
that these orbits are neither of finite dimension nor of finite codimension.

\proclaim{Theorem 4.3} a) For all $\ca\in\pg\Bbbk$ we have
$$
\oh^k_c(S_\nu,\ca) = 0 \ \ \text{if} \ \ k\neq 2\htt(\nu)\,.
$$ In particular, the functors $\oh^{2\htt(\nu)}_c(S_\nu,\ \,): \pg\Bbbk \to
\V$ are exact.

b) We have a natural equivalence of functors
$$
\Bbb H^* \ \cong \ \bigoplus_{\nu\in X_*(T)} \ \oh^{2\htt(\nu)}_c(S_\nu,\ )\ :
\
\pg\Bbbk
\to
\V
$$
\endproclaim

This result immediately gives the following consequence:

\proclaim{Corollary 4.4} The global cohomology functor $\Bbb H^*: \pg\Bbbk
\to
\V$ is exact.
\endproclaim

Here is a brief outline of the proof of theorem 4.3. Let us consider unipotent
radical
$\bar N$ of the borel $\bar B$ opposite to $B$. The $\bar N(\ck)$-orbits on
$\G$ are parametrized by $X_*(T)$: to each $\nu\in X_*(T)$ we associate the
orbit
$T_\nu = \bar N(\ck)\cdot \nu$\,. Recall that the $\GO$-orbits are
parametrized by $X_*(T)/W$. The orbits $S_\nu$ and $T_\nu$ intersect the
orbits $\G_\lambda$ as follows:
$$
\aligned &\text{a)} \ \ \dim(S_\nu\cap\G_\lambda) \ = \ \htt(\nu+\lambda)
\qquad
\text{if $\lambda$ is chosen dominant}
\\ &\text{b)} \ \ \dim(T_\nu\cap\G_\lambda) \ = \ -\htt(\nu+\lambda)
\qquad
\text{if $\lambda$ is chosen anti-dominant}
\\ &\text{c) \ \ the intersections in a) and b) are of pure dimension}\,.
\endaligned
\tag4.5
$$  In proving estimates  a) and b) we use the fact that the boundary $\partial
S_\nu$ is given by one equation in the closure $\bar S_\nu$. For the idea
behind the proof of c), see the appendix.  From the dimension estimates
(4.5a,b) above we conclude immediately that
$$
\aligned &\oh^k_c(S_\nu,\ca)\  =\ 0 \ \ \ \text{if} \ k>2\htt(\nu)
\\ &\oh^k_{T_\nu}(\G,\ca) \ =\ 0 \ \ \ \text{if} \ k<2\htt(\nu)\,.
\endaligned
\tag4.6
$$ Theorem 4.3 follows immediately from (4.6) and the following statement:
$$
\oh^k_c(S_\nu,\ca) \ = \oh^k_{T_\nu}(\G,\ca) \ \ \ \ \ \text{for all $k$}\,.
\tag4.7
$$  To see (4.7) we use the fact that $N(\ck)$-orbits and $\bar N(\ck)$-orbits
are in general position with respect to each other.

\remark{Remark 4.8} The decomposition of functors in theorem 4.3b is
independent of the choice of $N$. In the case of $N$ and its opposite
unipotent subgroup  $\bar N$  the corresponding decompositions are explicitly
related by $\oh^k_{S_\nu}(\G,\ca) \cong
\oh^k_{T_{w_0\cdot\nu}}(\G,\ca)$, where $w_0$ is the longest element in the
Weyl group. From this, and (4.7), we conclude that we could state theorem 4.3
replacing the functors $\oh^{2\htt(\nu)}_c(S_\nu,\ \,)$ by the equivalent set
of functors $\oh^{2\htt(\nu)}_{S_\nu}(\G,\ \,)$, where
$\oh^{2\htt(\nu)}_c(S_\nu,\ \,) \cong
\oh^{-2\htt(\nu)}_{S_{w_0\cdot\nu}}(\G,\
\,)$.
\endremark

\remark{Remark 4.9} The decomposition of $\G_\lambda$ into
$N(\ck)$-orbits and
$\bar N(\ck)$-orbits is an example of a perverse cell complex. Perverse cell
complexes are the analogues of CW-complexes for computing cohomology of
perverse sheaves instead of the ordinary cohomology. In the case at hand we
are in the situation analogous to the one for CW-complexes where  the
dimensions of all cells are of the same parity. We will develop the general
theory of perverse cell complexes elsewhere.
\endremark

\vskip 1cm

\subheading{\bf 5. The dual group}
\vskip .5cm

We will now apply Tannakian formalism as in \cite{DM} to $\pg\Bbbk$ and the
functor $\Bbb H^*$. In sections 2 and 3 we have given a tensor product
structure on the category $\pg\Bbbk$ via convolution and we have given
functorial associativity and commutativity constraints for this tensor
product.  To see that $\pg\Bbbk$ is a rigid tensor category, we still must
exhibit the identity object and construct duals. The identity object is given
by the sky scraper sheaf supported on the point
$1\cdot\GO\in\G$ whose stalk is
$\Bbbk$. The dual
$A\spcheck$ of a sheaf $A\in\pg \Bbbk$ is given as follows. Consider the
following sequence  of maps
$$
\G @<{\ \pi\ }<< \GK @>{\ i\ }>> \GK @>{\ \pi\ }>> \G\,,
\tag5.1
$$ where $i$ is the inversion on $\GK$, i.e., $i(g) = g^{-1}$. We define an
equivalence
$$
\iota:\pg \Bbbk \to \pg \Bbbk \text{ by} \ \iota(A) = \pi_* \tilde A \text{
where } i^*\tilde A = \pi^*A\,.
$$ Then the dual $A\spcheck$ is given by $A\spcheck = \iota (\Bbb D A)$, where
$\Bbb D$ denotes the Verdier dual.

In 4.1 we showed that the tensor product gets taken to the ordinary tensor
product in $\V$ by the functor $\Bbb H^*$. Furthermore, the associativity and
the commutativity constraints on $\pg\Bbbk$ get mapped to the standard ones on
$\V$ by $\Bbb H^*$. Corollary 4.4 says that $\Bbb H^*$ is exact and from this
it is not hard to deduce that it is also faithful. Thus, we have verified that
$\pg\Bbbk$ together with $\Bbb H^*$ constitutes a neutral Tannakian category
and by \cite{DM, theorem 2.11} we conclude:

\proclaim{Proposition 5.1} There exists an affine group scheme $\check G$ such
that the tensor category $\pg\Bbbk$ is equivalent to the (tensor category) of
representations of $\check G$. This equivalence is given via the fiber functor
$\Bbb H^*$.
\endproclaim

We claim:
\proclaim{Proposition 5.2} The affine groups scheme $\check G$ is isomorphic
to the Langlands dual of $G$.
\endproclaim

To see this, one may argue as follows. First of all, it is not difficult to
see that
$\check G$ is connected. By theorem 4.3b) we conclude that the dual torus
$\check T$ of $T$ is contained in $\check G$ and then one shows, as in
\cite{G}, that the torus $\check T$ is maximal. If $\check G$ had a unipotent
radical, then the category of $\pg\Bbbk$ would have certain non-trivial self
extensions of objects and this can easily be ruled out (this argument is due to
Soergel). As one can express the root datum of a reductive group in terms of
its irreducible representations one concludes, following \cite{G},  that
$\check G$ is the dual group of $G$.

A few remarks are in order. Because $\check G$ is reductive, one concludes
immediately that $\pg\Bbbk$ is semisimple. One can also see directly that
$\pg\Bbbk\cong \ps(\G,\Bbbk)$ is semisimple, for example from \cite{Lu,
theorem 11c}.

Let us make the statements of propositions 5.1 and 5.2 more concrete. Let
$\lambda\in X_*(T)/W = X^*(\check T)/W$. To $\lambda$ we can associate an
irreducible representation $V_\lambda$ of the Langlands dual group
$\check G$ on one hand, and a $\GO$-orbit $\G_\lambda$, and thus an
irreducible perverse sheaf $\Cal V_\lambda= j_{!*} \Bbbk_\lambda[\dim
\G_\lambda]$, $j:\G_\lambda\hookrightarrow \G$, on the other. Under the
equivalence of proposition 5.1 the sheaf $\Cal V_\lambda$ and the
representation
$V_\lambda$ correspond to each other. Furthermore,  the representation space
of $V_\lambda$ gets identified with the global cohomology of $\Cal V_\lambda$,
i.e., $V_\lambda = \Bbb H^*(\G,\Cal V_\lambda)$. This interpretation gives a
canonical basis for $\Cal V_\lambda$ as follows. From theorem 4.3, the fact
that $j_{!*}\Bbbk_\lambda[\dim
\G_\lambda] = \break ^p j_{!}\Bbbk_\lambda[\dim \G_\lambda]$, and (4.5c) we
conclude:
$$
\gathered
\Bbb H^k(\G,\Cal V_\lambda) \ = \  \bigoplus\Sb\nu\in X_*(T)\\k =
2\htt(\nu)\endSb
\oh^{2\htt(\nu)}_c(S_\nu,\Cal V_\lambda) \ = \\ \bigoplus\Sb\nu\in X_*(T)\\k =
2\htt(\nu)\endSb\oh^{2\htt(\lambda+\nu)}_c(S_\nu\cap\G_\lambda,\Bbbk)
\ = \ \bigoplus\Sb\nu\in X_*(T)\\k = 2\htt(\nu)\endSb
\Bbbk[\operatorname{Irr}(S_\nu\cap\G_\lambda)]\,.
\endgathered
\tag5.3
$$ Here $\Bbbk[\operatorname{Irr}(S_\nu\cap\G_\lambda)]$ denotes the vector
space spanned by the irreducible components of
$S_\nu\cap\G_\lambda$\,. Thus we get
$$ V_\lambda \ = \ \Bbb H^*(\G,\Cal V_\lambda) \ = \ \bigoplus\Sb\nu\in X_*(T)
\endSb \Bbbk[\operatorname{Irr}(S_\nu\cap\G_\lambda)]\,.
\tag5.4
$$ Note that the results above imply that the cohomology group $\Bbb
H^*(\G,\Cal V_\lambda)$ is generated by algebraic cycles.

\vskip 1cm

\subheading{\bf 6. Appendix}
\vskip .5cm

In this appendix we outline the proofs of theorem 2.7 and the statement
(4.5c).  Theorem 2.7 follows from the estimate:
$$
\dim[	m\inv S_\nu\ \cap\ p\inv(\GG_\la)  \bb{\GO}\to\tim \GG_\mu  ]
\le
\htt(\la+\mu+\nu)
\tag6.1
$$ for  coweights $\la,\mu, \nu\in X_*(T)$ such that
$\la$ and $\mu$ are dominant and
$\nu\in\barr{\GG_{\la+\mu} }$. The statement (6.1) can be proved exactly the
same way as the estimates (4.5a,b). We first directly verify (6.1) in the two
cases when $\nu=\lambda+\mu$  is dominant and  when $\nu= w_0(\lambda+\mu)$ is
antidominant ($w_0$ the longest element in the Weyl group). Then we use the
fact that the boundary $\partial S_\nu$ is given by one equation in  the
closure $\bar S_\nu$.

The  proof of the estimate (4.5c) is more involved as we use  a   Poisson
structure on the  ind-variety $\GG$. We choose an invariant non-degenerate
bilinear  form $\chi$ on $\fg$ and define an invariant non-degenerate form $(\
\,  , \ )$ on $\fg_\KK$ by the formula
 $(x,y)=\text{Res}\ \chi(x,y)$ for $ x,y\in\fg_\KK$.  The pair $\left(
\fg_\KK,\ (-,-)\right),$ has a Manin decomposition
 $(\fg_\KK)_+=\fg_\OO$ and $(\fg_\KK)_-= \fg_{z\inv\C[z\inv]}$, see, for
example, \cite{Dr}. This  formally defines
 a Poisson structure on the  ind-group $G(\KK)$ which descends to a Poisson
structure on $\GG= G(\KK)/\GO$. We have:

\proclaim{ Lemma 6.2} (a) The symplectic leaves in $\GG$ are the
intersections  of $\GO$-orbits  and the orbits of the negative congruence
subgroup $K_- =_{\text{def}} G({z\inv\C[z\inv]})$.

(b) The $N(\KK)$-orbits  are  coisotropic subvarieties of the Grassmannian
$\GG$.
\endproclaim

For a coweight $\nu\in X_*(T)$ we write $\GG^\nu=K_-\cd\nu\sub \GG$ for its
orbit under the negative congruence subgroup $K_-$. When the coweight $\nu$ is
antidominant the intersection $S_\nu\ \cap\ \GG_\la$ is a Lagrangian
subvariety of the symplectic leaf $ \GG^\nu\cap\GG_\la$. This implies
 (4.5c) in the antidominant case. To deduce (4.5c) for  general $\nu\in
X_*(T)$ we use the factorization
$$ S_\nu\ \cap\ \GG_\la\cong (S_\nu\ \cap\ \GG^\nu\cap\ \GG_\la)\tim (S_\nu
\cap\ \GG_\nu)
$$ and observe that the first factor
$S_\nu\ \cap\GG^\nu\ \cap\ \GG_\la$ is a Lagrangian  subvariety of the
symplectic leaf $ \GG^\nu\cap\GG_\la$.

\remark{Remark 6.4} Using the same techniques one can also prove a stronger
form of estimate (6.1). Namely, that the variety $m\inv S_\nu\cap
p\inv(\GG_\la)  \tim _\GO\GG_\mu $ is of pure dimension
$\htt(\la+\mu+\nu)$.
\endremark

\enddocument